\documentclass[12pt]{amsart}
\usepackage{amstext}
\usepackage{amssymb}
\usepackage{amsmath}
\usepackage{latexsym}
\usepackage{graphicx}
\usepackage{euscript}
\usepackage[latin1]{inputenc}
\usepackage{epsfig}
\usepackage{enumerate}

\def\centerbmp#1#2#3{\vskip#2\relax\centerline{\hbox to#1{\special
  {bmp:#3 x=#1, y=#2}\hfil}}}

\def\r{\mathbb{R}}

\def\c{\mathbb{C}}
\def\s{\mathbb{S}}

\def\ri{{\rm i}}
\def\re{{\rm e}}

\newtheorem{lemma}{Lemma}
\newtheorem{remark}{Remark}
\newtheorem{theorem}{Theorem}

\renewcommand{\qed}{\hfill $\square$}

\title[A Uniqueness theorem]{A Uniqueness Theorem for the singly periodic genus-one helicoid$^\ast$}
\author[A.\ Alarcón]{Antonio Alarcón }
\author[L.\ Ferrer]{Leonor Ferrer}
\author[F.\ Martín]{Francisco Martín}

\thanks{$\ast$ Research partially
supported by MCYT-FEDER grant number MTM2004-00160. \newline
 }

\address{Departamento de Geometría y Topología\hfill\break\indent  Universidad de Granada, \hfill\break\indent18071, Granada \hfill\break\indent Spain}

\email{alarcon@ugr.es, \,\,\,  lferrer@ugr.es, \,\,\, fmartin@ugr.es}

\begin{document}

\begin{abstract} The singly periodic genus-one helicoid was in the origin of the discovery of the first example of a complete minimal surface with finite topology but infinite total curvature, the celebrated Hoffman-Karcher-Wei's genus one helicoid. The objective of this paper is to give a uniqueness theorem for the singly periodic genus-one helicoid provided the existence of one symmetry.
 \end{abstract}

\maketitle
\section{Introduction} \label{sec:intro}
In the last few years, one of the most active focus in the study of minimal surfaces has been the genus-one helicoid. The existence of such a surface was proved by D. Hoffman, H. Karcher and F. Wei in \cite{hkwb} and it was at that moment the first example of an infinite total curvature but finite topology embedded minimal surface.

One important step in the discovery  of the genus-one helicoid was the construction by D. Hoffman, H. Karcher and F. Wei in \cite{hkw} of a singly periodic minimal surface which is invariant under a translation so that the quotient has genus one. This minimal surface is  called the singly periodic genus-one helicoid and  it will be represented as ${\mathcal H}_1$. Other that the helicoid itself, this example was the first embedded minimal surface ever found that in asymptotic to the helicoid. The helicoid ${\mathcal H}_1$ belongs to a continuous family of twisted periodic helicoids with handles that converges to a genus one helicoid. The continuity of this family of surfaces and the subsequent embeddedness of this genus one helicoid were obtained by D. Hoffman, M. Weber and M. Wolf in \cite{hww, weber}. Although there are numerical evidences that there is only one embedded helicoid with genus one, to our knowledge this fact remains unproven.
 \begin{figure}[ht] 
	\begin{center}
	\scalebox{0.35}{\includegraphics{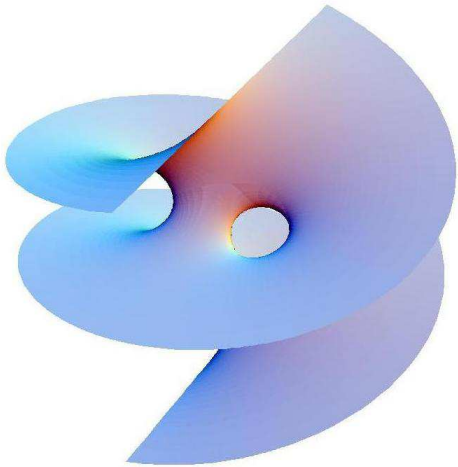}}
		\end{center}
\caption{Helicoid  ${\mathcal H}_{1}$.}\label{g1}
\end{figure}

Furthermore, a recent result by W.H. Meeks and H. Rosenberg asserts that any properly immersed minimal surface with finite topology and one end must be asymptotic to a helicoid with handles and can be described by its Weierstrass data $(\frac{dg}{g},dh)$ on a compact Riemann surface (see \cite{mra}).

From the preceding arguments, uniqueness results about ${\mathcal H}_1$ and the family of twisted periodic helicoids with handles derived from it become more interesting. In this setting, L. Ferrer and F. Martín have obtained recently the following result:
\begin{theorem} \label{hkw} Any complete, periodic, minimal surface containing a vertical line, whose quotient by vertical translations has genus one, contains two parallel horizontal lines, has two helicoidal ends and total curvature $-8 \pi$ is ${\mathcal H}_1$.
\end{theorem}
This result was essentially obtained in \cite[Theorem 1]{hkw} where D. Hoffman, H. Karcher and F. Wei proved that a surface with the qualitative properties of the surface in Theorem \ref{hkw} belongs to a two-parameter family of Weierstrass data and the  period problem is solvable in this family. Our contribution consists of giving a new approach to the proof of the uniqueness of the period problem (see Remark 3 in \cite{fm2}). Indeed, we prove that there is only one pair of this parameters that solves the period problem.

The main objective of the present paper is to prove the following uniqueness theorem for ${\mathcal H}_1$ that improves the aforementioned one. 
\begin{theorem} \label{teorema} Any properly embedded, singly periodic minimal surface that is symmetric respect to a vertical line, whose quotient by a vertical translation  has genus one,  two helicoidal ends and total curvature $-8 \pi$ is ${\mathcal H}_1$.
\end{theorem}
In order to demonstrate this result we will see that a surface satisfying the hypothesis of Theorem \ref{teorema} also verifies the hypothesis of Theorem \ref{hkw} and so our result is a direct consequence of the previous one.

\section{Preliminaries} 
Given $X=(X_1,X_2,X_3):M \longrightarrow \r^3$ a conformal minimal immersion we denote by $g:M \longrightarrow \overline{\c}=\c \cup \{\infty \}$ its stereographically projected Gauss map that is a meromorphic function and by $\Phi_3$ the holomorphic differential defined as $\Phi_3=d X_3+ {\rm i}d X_3^*$, where $X_3^*$ denotes the harmonic conjugate function of $X_3$. The pair $(g,\Phi_3)$ is usually referred to as the Weierstrass data of the minimal surface, and the minimal immersion $X$ can be expressed, up to translations, solely in terms of these data as 
\begin{equation}\label{eq:inmersion} \hspace{-0.5in} X={\rm Re}\int^z (\Phi_1,\Phi_2,\Phi_3)={\rm Re}\int^z \left(\frac{1}{2} \left( \frac{1}{g}-g \right),\frac{{\rm i}}{2}\left(\frac{1}{g}+g \right),1 \right) \Phi_3 \;, \end{equation}
where ${\rm Re}$ stands for real part and $z$ is a conformal parameter on $M$. The pair $(g,\Phi_3)$ satisfies certain compatibility conditions:
\begin{enumerate}
\item[{ i)}] The zeros of  $\Phi_3$ coincide with the poles and zeros of  $g$, with the same order.
\item[{ ii)}] For any closed curve $\gamma \subset M$,
\begin{equation} \label{ec:periodos} \overline{\int_\gamma g \Phi_3}=\int_\gamma \frac{\Phi_3}{g} \quad, \quad {\rm Re} \int_\gamma \Phi_3=0 \; . \end{equation}
\end{enumerate}
Conversely, if $M$ is a Riemann surface, $g:M \to \overline{\c}$ is a meromorphic function and $\Phi_3$ is a holomorphic one-form on $M$ fulfilling the conditions i) and ii), then the map $X:M \to \r^3$ given by
\eqref{eq:inmersion} is a conformal minimal immersion with Weierstrass data $(g,\Phi_3)$.

Condition ii) stated above deals with the independence of  \eqref{eq:inmersion} on the integration path, and it is usually  called the period problem.

\section{Proof of Theorem \ref{teorema}}\label{sec:3}

Let $\widetilde{X}:\widetilde{M} \to \r^3$ be a minimal immersion satisfying the conditions in Theorem \ref{teorema}. We label $\mathsf{t}$ as the vertical translation and $R$ as the axis of symmetry. Denote by $M=\widetilde{M}/\langle \tau \rangle$ and by $X:M \to \r^3/\langle \mathsf{t} \rangle$ the immersion that verifies $X \circ p= \widetilde{p} \circ \widetilde{X}$, where $\tau$ is the biholomorphism in $\widetilde{M}$ induced by $\mathsf{t}$, and $p:\widetilde{M} \to M$ and $\widetilde{p}:\r^3 \to \r^3/\langle \mathsf{t} \rangle$ are the canonical projections. Note that the Weierstrass data of the immersion $\widetilde{X}$, that we denote $(g,\Phi_3)$, can be induced in the quotient $M$. We also denote $(g,\Phi_3)$ as the induced Weierstrass data.

From general results (\cite{mrc}) and our assumptions we know that  $M$ is conformally equivalent to a torus $T$ minus two points $\{E_1,E_2\}$ and the Weierstrass data extend meromorphically to the ends $\{E_1,E_2\}$. 

 Since we are assuming that $\{E_1,E_2\}$ are helicoidal ends with vertical normal vectors we have that $\Phi_3$ has simple poles with imaginary residues at these points. Up to an homotethy we can assume that  $\text{Res}(\Phi_3 ,E_1)=-\text{Res}(\Phi_3 ,E_2)=\ri$. Furthermore, as $T$ is a torus,  there exist two zeros  $\{V_1,V_2\}$  of $\Phi_3$ in $T$ and thereby the divisor of $\Phi_3$ is given by 
\begin{equation} \label{divisorfi}(\Phi_3)=\frac{V_1 V_2}{E_1 E_2} \;. \end{equation}

Concerning $g$, using the formula 
$$ \int_{M} K dA_e=-4 \: \pi \, \deg(g) \; ,$$
we obtain that $\deg(g)=2$. Hence $g$ has two zeros and two poles that must coincide with the points $\{V_1,V_2,E_1,E_2\}$. Since  the normal vectors at the ends have opposite directions, up to relabeling, we can assume that the divisor of $g$ is
\begin{equation} \label{divisorg} (g)=\frac{V_1 E_1}{V_2 E_2} \;. \end{equation}
Moreover $V_1 \neq V_2$. If not $\deg(g)<2$. 

We shall call ${\mathcal S}_3$ the isometry of $\widetilde{M}$ induced by the symmetry of 180º rotation about the line $R$ and $S_3:T \to T$ the involution induced by ${\mathcal S}_3$ on the torus.

Now we study the intersection of $\widetilde{X}(\widetilde{M})$ with the horizontal planes.
\begin{lemma}\label{lem:alpha} 
 $\alpha_k=\widetilde{X}(\widetilde{M}) \cap \{x_3=k\}$ is either a simple curve $\ell_k$ diverging to both ends  and containing the point $p_k=R \cap \{x_3=k\}$ or the union of such a curve $\ell_k$ and a Jordan curve that cuts $\ell_k$ orthogonally at two points. Moreover, this situation occurs once in each fundamental piece. As a consequence we obtain $R \subset \widetilde{X}(\widetilde{M})$.
\end{lemma}
\proof Firstly we recall that the intersection of a minimal surface with a plane is, in a neighborhood of each point $p$, a set of $n$ analytic curves that intersect each other at $p$ at
an angle $\pi/n$. Furthermore, if the plane is the tangent plane to the minimal surface at $p$ then, the multiplicity
of the Gauss map of the surface at $p$ is $n-1$ if, and only if,  the plane intersects the surface along $n$ curves  in
a neighborhood of $p$. Obviously, the multiplicity of the Gauss map at $p$ is $1$, if and only if, the tangent plane to a minimal surface  intersects the surface along two orthogonal curves.

\begin{figure}[ht]
\begin{center}
\scalebox{0.6}{\includegraphics{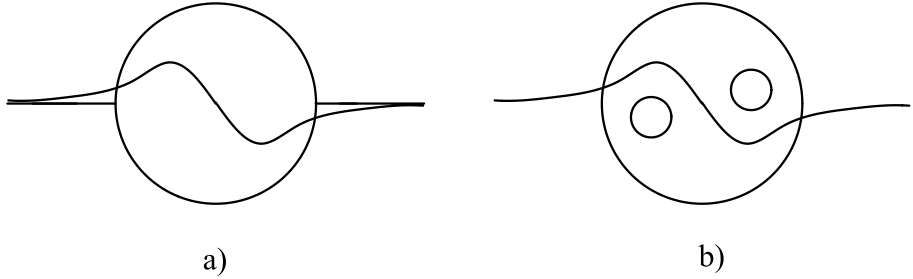}} 
\end{center}
\caption{The curve $\alpha_k$.}\label{circulos}
\end{figure}

We restrict ourselves to the fundamental piece $M$. Since there are only two points with vertical normal vector in $M$,  we deduce that the intersection with any horizontal plane consists of a set of disjoint simple analytic curves except for the plane $\{x_3=k_0\}$ that contains the points $V_1$ and $V_2$. Observe that these points must be at the same horizontal plane by the symmetry. Moreover, as the ends are of helicoidal type we have that, outside a sufficiently large vertical cylinder, the curve $\alpha_k$ has two symmetric connected components that diverge to both ends.  Observe that for $k\neq k_0$ there exists a simple curve $\ell_k \subset \alpha_k$ that contains these two connected components  (see Fig. \ref{circulos}.a).

Moreover, for any $k\neq k_0$ we can deduce that the point $p_k=R \cap \{x_3=k\} \in \ell_k$. Otherwise, the point $p_k$ would be contained in one of the two connected components of $\{x_3=k\} \setminus \ell_k$. Note that outside the cylinder the symmetry interchanges those connected components. As $p_k$ and $\ell_k$ are invariant under the symmetry of 180º rotation about the line $R$, we deduce that $p_k$ also belongs to the other connected component, which is a contradiction. Since  $R \cap \widetilde{X}(\widetilde{M})$ is an continuous set, we deduce that  $S_3$ is an antiholomorphic involution. Hence the set of fixed points is a whole curve and so $R \subset \widetilde{X}(\widetilde{M})$.

In relation with the behavior of $\alpha_k$ in the interior of the cylinder, firstly we shall prove that $\alpha_{k}$ does not contain bounded connected components.

\begin{figure}[ht]
\begin{center}
\epsfig{file=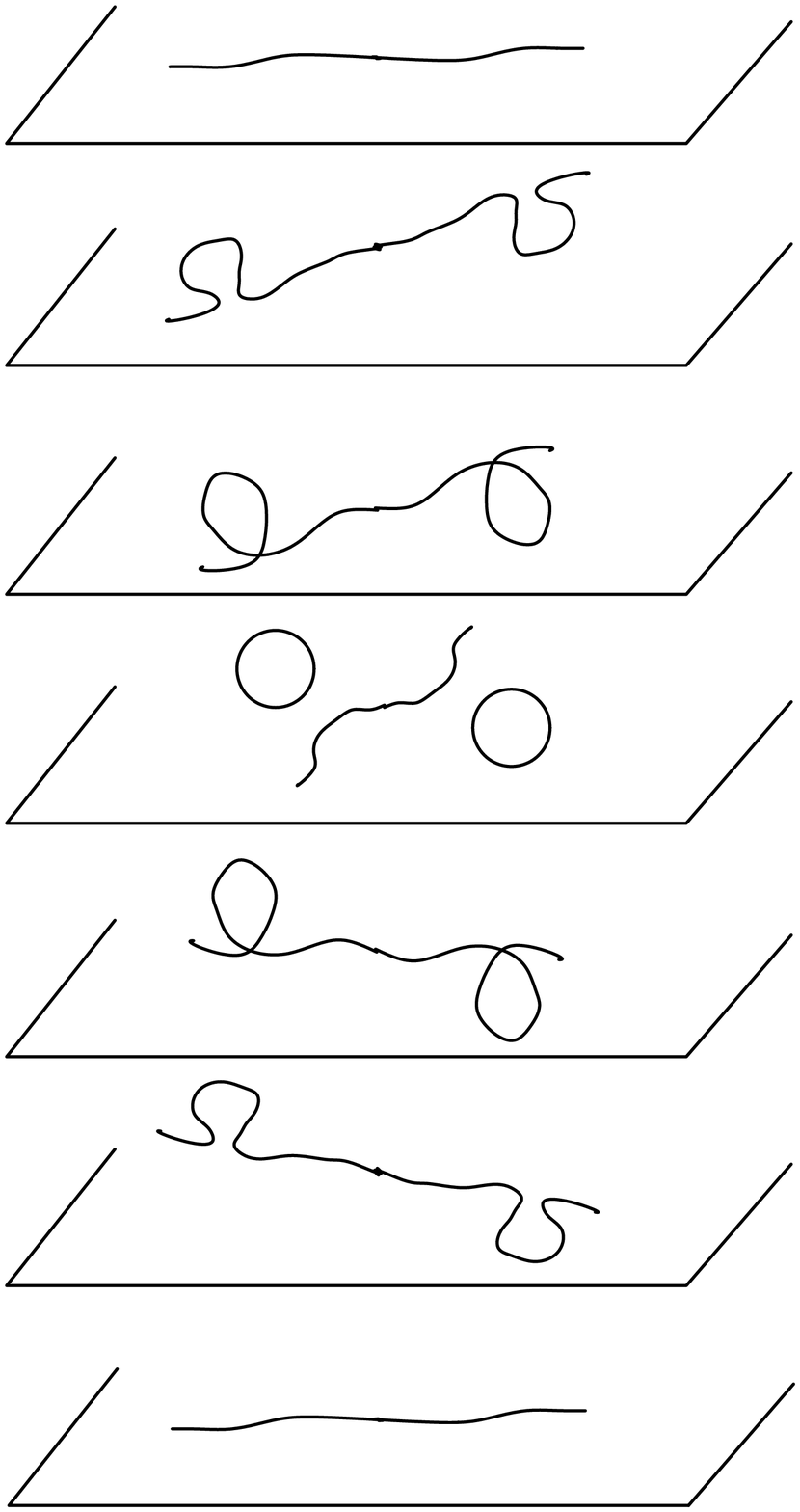,height=10cm,width=10cm}
\end{center}
\caption{}\label{fig:evolucion}
\end{figure}

Assume that there exists a bounded connected component. Then, by the symmetry, there exist at least two of these connected components (see Fig. \ref{circulos}.b). If these curves would appear in  $\alpha_k$, for any $k$, then  $\widetilde{X}(\widetilde{M})$ does not divide $\r^3$ in two connected components, contradicting the embeddedness of the surface. Consequently, for some $k_1$ the curve $\alpha_{k_1}$ must be connected. Taking into account the symmetry and the first paragraph in this proof we infer that the evolution of the curves $\alpha_k$ is as in Fig. \ref{fig:evolucion}. But this contradicts the fact that there are only two points of vertical normal vector in $M$.

Therefore, for $k\neq k_0$ the intersection $\alpha_k=\ell_k$, it is to say, it is a simple curve  symmetric respect to $R$ and diverging to both ends. 

\begin{figure}[ht]
\begin{center}
\epsfig{file=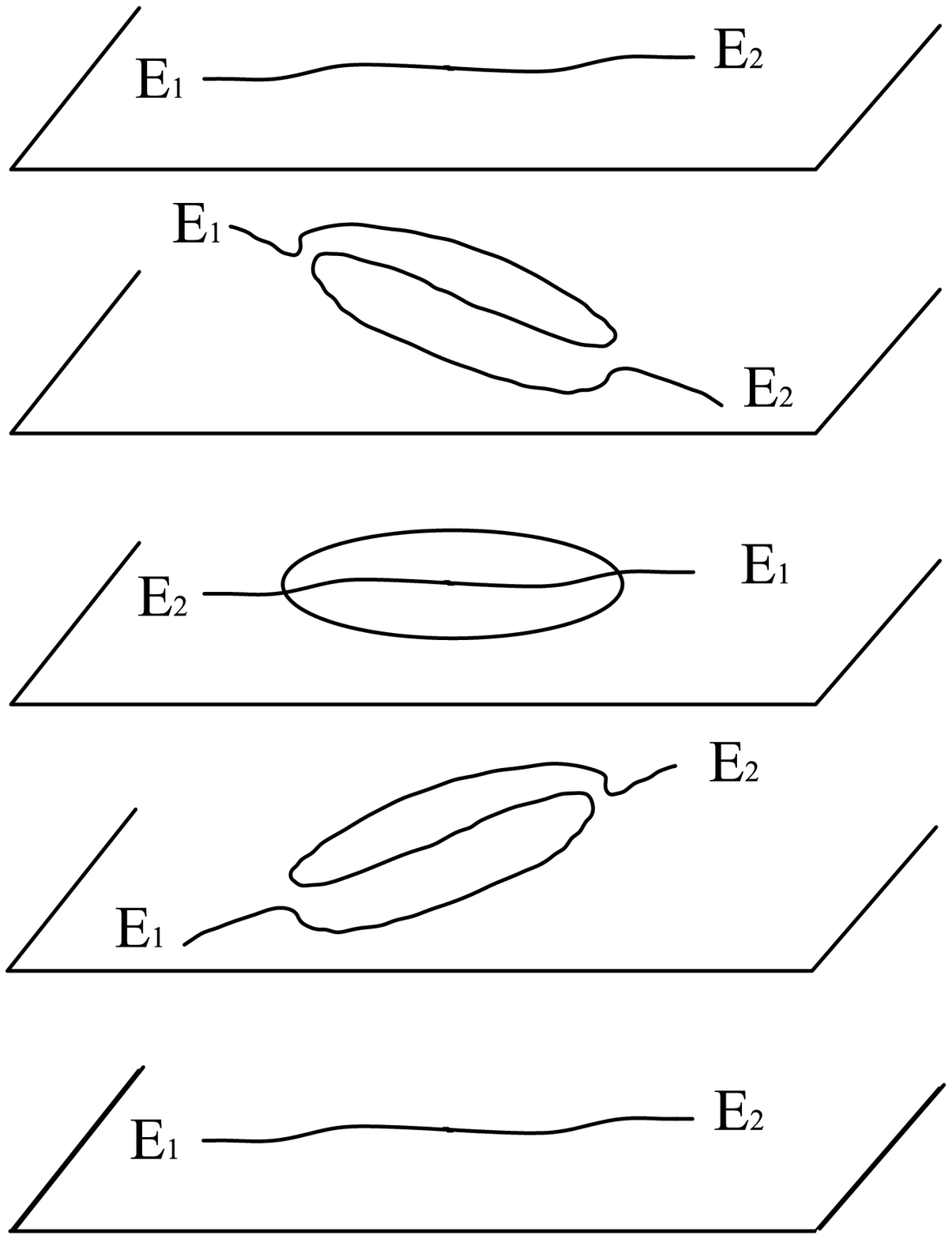,height=11cm,width=10cm}
\end{center}
\caption{}\label{fig:evoluciona}
\end{figure}
On the other hand the curve $\alpha_{k_0}$ must consist of two curves: $\ell_{k_0}$ as the previous ones and another Jordan curve $\beta$ that intersects $\ell_{k_0}$ orthogonally at the points $V_1$ and $V_2$. Indeed, the evolution of the curves $\alpha_k$ is as represented in Fig. \ref{fig:evoluciona}. \qed

\vspace{0.12in}

Henceforth, we shall assume that the vertical line contained in $\widetilde{X}(\widetilde{M})$ is the axis $\{x=y=0\}$ and we shall call $L$ the closed curve $L \subset M$ such that $X(L)= \{x=y=0\}\cap \r^3/\langle \mathsf{t} \rangle$.

\begin{remark} In general, the horizontal level curves of a minimal annular end that is asymptotic to a helicoid are not asymptotic to straight lines. One example of this situation can be found in \cite{hoffman}, Remark 4.

However, in our case it is easy to see that the curves $\ell_k$ converge to straight lines. Indeed, it is known that if $E$ is an end of a properly embedded, singly periodic minimal surface of genus one and invariant under a vertical translation and we assume $g$ has a zero at the end, it is possible to consider a conformal coordinate $z$ around the end such that
$$ g(z)=z h(z) \; , \Phi_3(z)=\frac{-i}{z} dz \; ,$$
where $h$ is a holomorphic function at the end with $h(0)=a_0 \neq 0$. 
Hence, we obtain that the projection of the end over the $(x_1,x_2)$-plane is given by
$$ (x_1+\ri x_2)(z)=c+\frac{-\ri }{2 \overline{a_0 z}}+O(z) \; ,$$
where $c \in \c$. We recall that $(\rm{Re}(c),\rm{Im}(c),0)+\langle (0,0,1)\rangle$ is the axis of the helicoidal end. As we are assuming that this axis is the $x_3$-axis we have that $c=0$.  For more details see \cite{joaquin}.

Taking into account that $ x_3(z)=  -\ri \log(z)$ and
the above expression, we obtain that the curve $\ell_k$ can be parametrized in a neighborhood of the end as
$$ \ell_k(r)=\frac{1}{2}\left(\frac{  1}{a r}\sin(k+\theta_0)+O(r),\frac{  -1}{a r}\cos(k+\theta_0) +O(r)\right) \; ,$$
where $r\in ]0,\varepsilon[$ and $a_0=a \re^{\theta_0 \ri}$. Clearly, this curve is asymptotic to the line 
$$\frac{1}{2} \left( \frac{  1}{a r}\sin(k+\theta_0), \frac{  -1}{a r}\cos(k+\theta_0 )\right) \;.$$
From the above argument, Fig. \ref{circulos}.a) is a realistic representation of the curves $\ell_k$. 
\end{remark}
\vspace{0.12in}

In the proof of Lemma \ref{lem:alpha} we obtained that $S_3:T \to T$ is an antiholomorphic involution. Then, it is not hard to see that the symmetry acts on the one-forms $\Phi_i$ as follows
\begin{equation} \label{ec:phis} S_3^\ast(\Phi_i)=-\overline{\Phi_i} \;, i=1,2, \quad  S_3^\ast(\Phi_3)=\overline{\Phi_3} \;.
\end{equation}
Hence taking into account that $g=\frac{\Phi_3}{\Phi_1-\ri \Phi_2}$ we obtain 
\begin{equation}\label{ec:gs} g(S_3(p))=\frac{1}{\overline{g(p)}} \; , p \in T \;. \end{equation}

\begin{lemma} \label{lem:g} The Gauss map $g$ has exactly two ramification points in $L$.
\end{lemma}
\proof Consider the set 
$${\mathcal G}=\{p \in M \mid |g(p)|=1\} \; ,$$
it is to say, the set of all the points in $M$ with horizontal normal vector. Since $\mathcal G$ is the nodal set of the harmonic function $\log(|g|)$ we have that it consists of a set of analytic curves. Moreover, when the nodal lines meet they form an equiangular system and the intersection points coincide with the ramification points of $g$.  Observe also that $\mathcal G$ is compact in $T$ and the curves in $\mathcal G$ do not converge to the ends. Thus $\mathcal G$ is compact in $M$.

Clearly $L \subset {\mathcal G}$. Note that Lemma \ref{lem:alpha} guarantees the existence of two points with horizontal normal vector in $\beta \setminus \{V_1,V_2\}$ and thereby $L \neq {\mathcal G}$.

\begin{figure}[ht]
\begin{center}
\scalebox{0.5}{\includegraphics{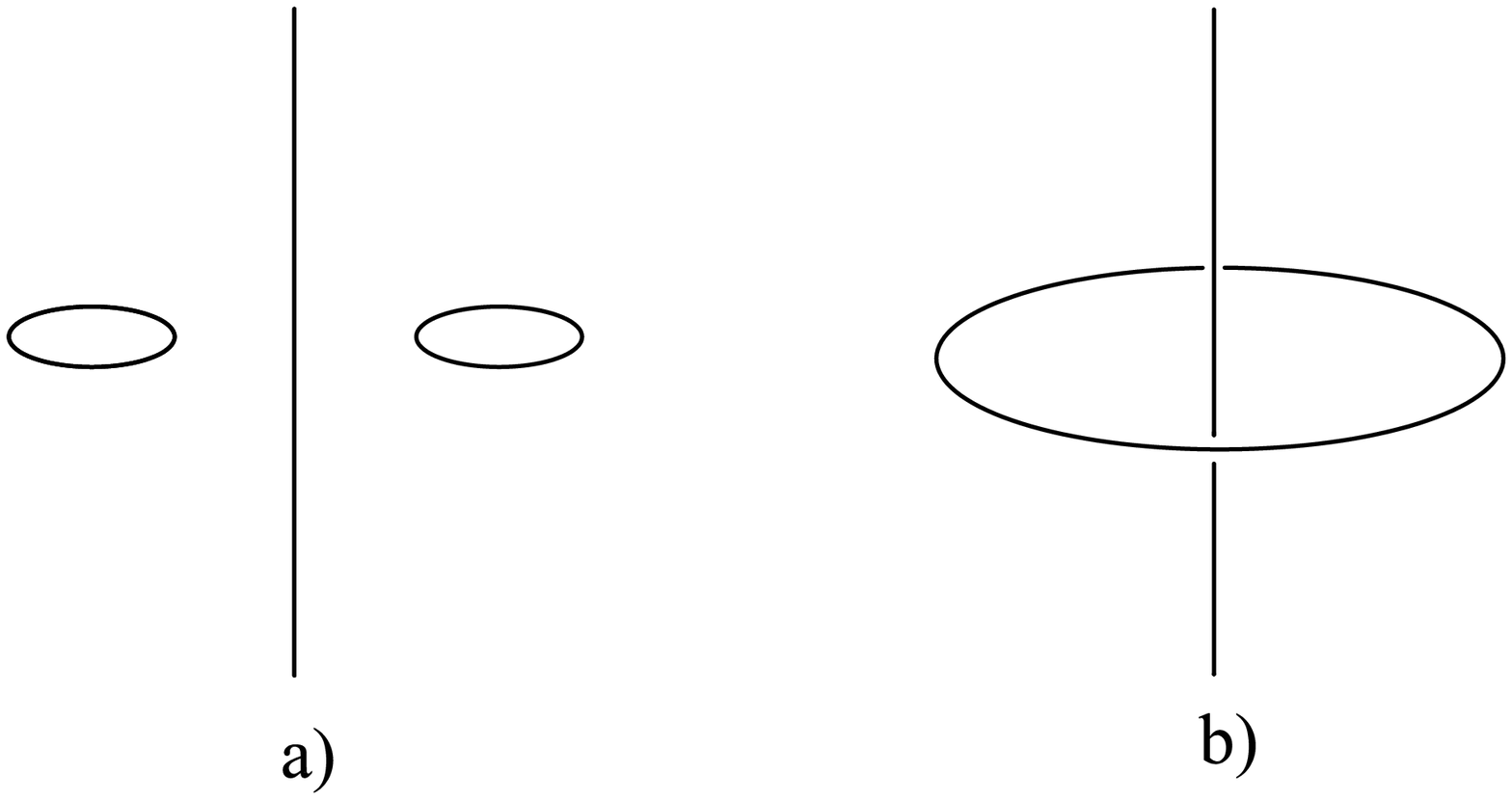}} 
\end{center}
\caption{}\label{piruli}
\end{figure}

Suppose that there are not ramification points of $g$ in $L$ and let $\gamma$ be a closed curve  in ${\mathcal G}$ different from $L$. From our assumptions we have that $\gamma \subset  {\mathcal G} \setminus L$ and $g:L \to \s^1$ is a bijection. Moreover, by the symmetry we have that $S_3(\gamma) \subset {\mathcal G}$. If $\gamma \neq S_3(\gamma)$, then the image of $\gamma$ and  $S_3(\gamma)$ by $g$ would cover twice  $\s^1 \subset \c$ contradicting $deg(g)=2$  (see Fig. \ref{piruli}.a). Suppose now that $\gamma = S_3(\gamma)$ (see Fig. \ref{piruli}.b).

Reasoning as before we obtain that $g:\gamma \to \s^1$ is a bijection. On the other hand, taking into account \eqref{ec:gs} we have
$$ g(S_3(p))=\frac{1}{\overline{g(p)}} \; , \quad p \in \gamma \; . $$
But $g(p) \in \s^1$ and so $g(S_3(p))=g(p)$ which contradicts the injectivity of $g|_\gamma$. 

\begin{figure}[ht]
\begin{center}
\scalebox{0.5}{\includegraphics{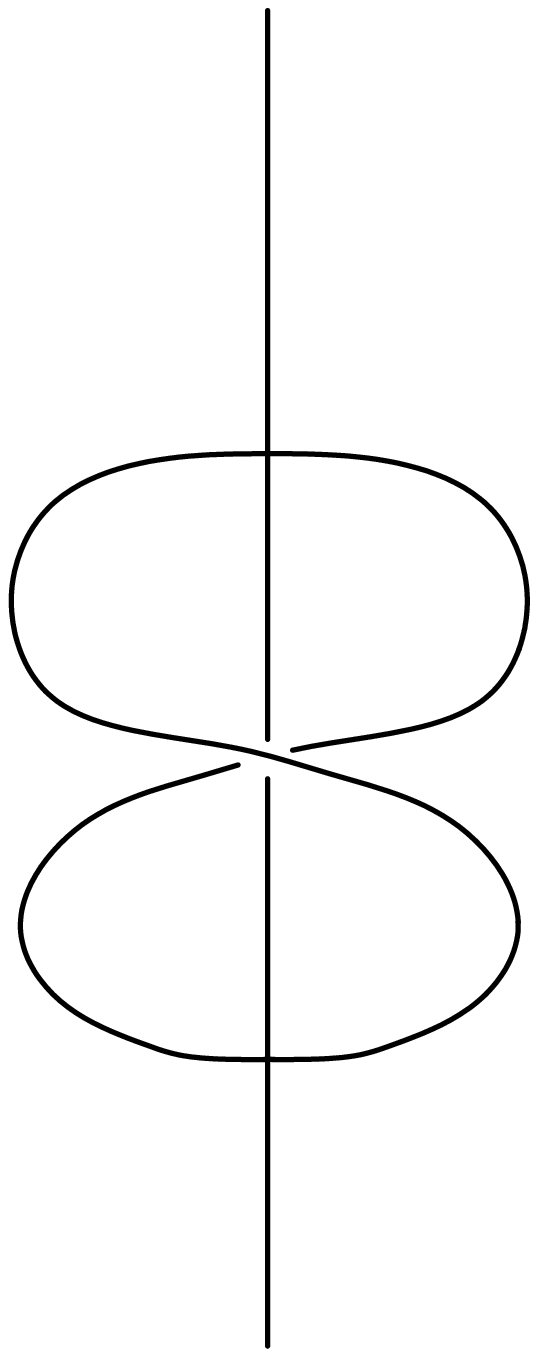}} 
\end{center}
\caption{}\label{buena}
\end{figure}

Then, there exist ramification points of $g$ in $L$.  Using that $\deg(g)=2$ and the Riemann-Hurwitz formula we deduce that $g$ has four ramification points in $T$ with  branch number one. Moreover, there are exactly two of these points in $L$ (see Fig. \ref{buena}) because any other situation leads us to a contradiction. \qed

Now we continue with the proof of the theorem. Hereafter we denote by $z=g:T \to \overline{\c}$. Since $z$ has four ramification points in $T$ with branch number one, it is known that the torus $T$ is conformally equivalent to the torus
$$\left\{ (z,w) \in \overline{\c}^2 \big| \: w^2=\frac{(z-a) (z-b)}{(z-c) (z-d)}\right\} \; .$$
where  $a,b,c,d \in \c$ are the images by $z$ of the ramification points of $g$. We identify $T$ with this torus. From Lemma \ref{lem:g} we have that two of these points are in $\s^1 \subset \overline{\c}$. Up to a rotation, we can assume that $c=\re^{\ri \rho}$, $d=\re^{-\ri \rho}$ and $-1 \in g(L)$.  Since $S_3$ preserves the set of ramification points of $g$, using \eqref{ec:gs} we have that $b=\frac{1}{\overline{a}}$. After these considerations $w$ can be written as
$$w^2=\frac{(z-a) (\overline{a}z-1)}{(z-\re^{\ri \rho}) (z-\re^{-\ri \rho})} \; ,$$
with $a \in \c$, $a \neq \re^{\ri \rho}$ and $a \neq \re^{-\ri \rho}$.

From \eqref{ec:gs} we have 
$$ S_3(z,w)=(\frac{1}{\overline{z}},\pm \overline{w}) \; .$$
Next we try to determine the appropriate sign of the second component. In order to do this, we observe that $S_3$ fixes the point 
$$\left( -1,+ \sqrt{\frac{1+|a|^2-2 \text{Re}(a)}{2-2 \cos(\rho)}}\right) $$ 
and  therefore the correct expression for $S_3$ is
$$ S_3(z,w)=(\frac{1}{\overline{z}},\overline{w}) \; .$$
Hence, up to relabeling, we have
$$ E_1=(0,\sqrt{a}) \; , V_1=(0,-\sqrt{a})\; , E_2=(\infty, \sqrt{\overline{a}})\; , V_2=(\infty, -\sqrt{\overline{a}})\; . $$

Next we prove that $a \in \r$.
\begin{lemma} If $M$ satisfies the hypothesis of Theorem \ref{teorema} then $a \in \r$.
\end{lemma}
\proof  Next we consider the meromorphic functions $ w-\sqrt{a}, w+\sqrt{a},w+\sqrt{\overline{a}}$ with divisors 
$$ (w-\sqrt{a})=\frac{E_1 P_0}{e^{i\rho} e^{-i\rho}} \; , \quad (w+\sqrt{a})=\frac{V_1 P_1}{e^{i\rho} e^{-i\rho}} \; , \quad (w+\sqrt{\overline{a}})=\frac{V_2 P_2}{e^{i\rho} e^{-i\rho}} \; ,$$ 
where $P_0=\big(z_1\,,\,+\sqrt{a}\big)\; , P_1=\big(z_1\,,\,-\sqrt{a}\big)\; ,
P_2=\left( \frac{1}{\overline{{z}_1}}\,,\,-\sqrt{\overline{a}}\right)$ and 
$$ z_1 = \frac{1+|a|^2-2a\cos\rho}{2\: \text{Im}(a)}\,i \; .$$
Note that $z_1 \neq 0$ because $a \neq \re^{\ri \rho}$ and $a \neq \re^{-\ri \rho}$. We introduce now the following meromorphic 1-form
$$\eta=\text{Im}(\sqrt{a})(z-z_1)\frac{w+\sqrt{\overline{a}}}{w-\sqrt{a}}\; \tau\; ,$$
where $\tau=\frac{dz}{(z-\re^{\ri \rho}) (z-\re^{-\ri \rho})w}$ is the holomorphic 1-form on the torus. It is easy to check that the divisor of $\eta$ is given by
$$ (\eta)= \frac{P_1 P_2}{E_1 E_2} \; ,$$
and $\text{Res}(\eta ,E_1)=\ri$.  Now, observe that $\Phi_3$ and $\eta$ are two meromorphic 1-forms on the torus with the same poles and the same residues at these poles. Consequently, the difference between them is a multiple of $\tau$, it is to say
\begin{equation} \label{suma} \Phi_3= \eta +\lambda \tau \; , \end{equation}
with $\lambda \in \c$. It is easy to see that 
\begin{equation} \label{etatau} S_3^\ast(\eta)=\overline{\eta}\; ,\quad  S_3^\ast(\tau)=-\overline{\tau} \; . \end{equation}
Recall that $S_3^\ast(\Phi_3)=\overline{\Phi_3}$. Then from \eqref{etatau} and \eqref{suma} we deduce that 
$$ \overline{\eta} -\lambda \overline{\tau}= \overline{\eta} + \overline{\lambda} \overline{\tau} \; .$$
Hence we obtain that $\text{Re}(\lambda)=0$. Thereby we write $\lambda= \ri r$ with $r \in \r$. Now we use that $\Phi_3$ has a zero at $V_1$. Substituting in the expression of $\Phi_3$ we get
$$0=  \eta(V_1) + \ri r \: \tau(V_1)=\left(\text{Im}(\sqrt{a})z_1 \frac{\sqrt{\overline{a}}-\sqrt{a}}{2  \sqrt{a}} + \ri r \right)\tau(V_1)\; .$$
Taking into account that $\tau(V_1)\neq 0$ and simplifying in the above equality we obtain
$$\text{Im}(\sqrt{a})z_1 \frac{\sqrt{\overline{a}}-\sqrt{a}}{2  \sqrt{a}} + \ri r=0 $$
which is equivalent to 
$$z_1 \frac{\sqrt{\overline{a}}-\sqrt{a}}{2  \sqrt{a}}= \overline{z_1} \frac{\sqrt{a}-\sqrt{\overline{a}}}{2 \sqrt{\overline{a}}}=-z_1 \frac{\sqrt{a}-\sqrt{\overline{a}}}{2 \sqrt{\overline{a}}} \; .$$
As $z_1 \neq 0$, from the above expression we get $\text{Im}(a)=0$ and we obtain then  $a \in \r$. \qed

Summarizing we have the following Weierstrass data 
$$  g=z \; , \quad   \Phi_3=\eta= -i \, \frac{1+a^2-2a\cos\rho}{2\sqrt{a}}\, \frac{w+\sqrt{a}}{w-\sqrt{a}}\; \tau \; , $$
on the torus $\left\{ (z,w) \in \overline{\c}^2 \big| \: w^2=\frac{(z-a) (a z-1)}{(z-\re^{\ri \rho}) (z-\re^{-\ri \rho})} \right\}$ punctured at the points $ E_1=(0,\sqrt{a})$ and $E_2=(\infty, \sqrt{a})$. 

Therefore, if we denote by $S:T \to T$ the symmetry given by $S(z,w)=(\overline{z},\overline{w})$, from the above expressions it is easy to check that
$$ g \circ S=\overline{g} \; , \quad S^\ast(\Phi_3)=-\overline{\Phi_3} \; .$$
Since $S$ is an antiholomorphic involution and the associated isometry is a reflection respect to a horizontal line, we deduce that this horizontal line lies on the surface. Consequently, we have all the hypothesis of Theorem \ref{hkw} and so we can conclude that $\widetilde{M}$ is the helicoid ${\mathcal H}_1$.

\end{document}